\documentclass[10pt,twoside]{article}
\usepackage{amssymb}
\usepackage[mathscr]{eucal}
\usepackage{eufrak}
\usepackage{amsmath}
\usepackage{mathrsfs}
\usepackage[colorlinks=true,
   urlcolor=blue,           filecolor=green,      
   citecolor=green,      
   linkcolor=red,           bookmarks=true,
  unicode,
   plainpages=false,   ]{hyperref}
\usepackage{color}

\def\n{\nabla}
\def\intl#1{\int\limits_{#1}}
\def\intll#1#2{\int\limits_{#1}^{#2}}
\def\dm{|\hskip-0.05cm|}
\def\OO{\Omega}
\def\displ{\displaystyle}
\def\VSE{\vspace{6pt}\\&\displ }
\def\VS{\vspace{6pt}\\\displ }
\def\rf#1{{\rm(\ref{#1})}}
\def\chiu{\hfill$\displaystyle\vspace{4pt}
\underset\Box\null$\par}
\def\R{\Bbb R}
\def\N{\Bbb N}

\def\à{\`{a}}

\def\vep{\varepsilon}

\def\be{\begin{equation}}
\def\ba{\begin{array}}
\def\ea{\end{array}}
\def\ee{\end{equation}}
\def\vs1{\vspace{1ex}}

\def\ov{\overline}
\def\po{\partial\Omega}
\def\Pr{{\bf Proof. }}
\def\é{\'{e}}
\font\sc=cmcsc10
\title{\Large A  note on the Prodi-Serrin conditions for the regularity  of a   weak solution to the Navier-Stokes equations}
\author{\sc   Paolo Maremonti
\thanks{
Dipartimento di Matematica e Fisica,  
Universit\`{a} degli Studi della Campania
``L. Vanvitelli'', via Vivaldi 43, 81100 Caserta,
 Italy.
paolo.maremonti@unicampania.it}}
\date{}
\begin{document}
\markboth{\footnotesize\rm P.
Maremonti} {\footnotesize\rm
A note on the Prodi-Serrin conditions...}
\maketitle \noindent{\bf Abstract} - {\small The paper is concerned with the regularity of weak solutions to the Navier-Stokes equations. The aim  is to investigate on a  relaxed Prodi-Serrin condition in order to obtain regularity for $t>0$. The most interesting aspect of the result is    that no  compatibility condition is required  to the initial data $v_\circ\in J^2(\OO)$.} 
 \par\noindent
 \vskip -0.7true cm\noindent
\newcommand{\red}{\protect\bf}
\renewcommand\refname{\centerline
{\red {\normalsize \bf References}}}
\newtheorem{ass}
{\bf Assumption} 
\newtheorem{defi}
{\bf Definition} 
\newtheorem{tho}
{\bf Theorem} 
\newtheorem{rem}
{\sc Remark} 
\newtheorem{lemma}
{\bf Lemma} 
\newtheorem{coro}
{\bf Corollary} 
\newtheorem{prop}
{\bf Proposition} 
\renewcommand{\theequation}{\arabic{equation}}
\setcounter{section}{0}
\section{Introduction}\label{intro}We consider the 3-d Navier-Stokes initial boundary value problem:
\be\label{NS}\ba{l}v_t+v\cdot
\nabla v+\nabla\pi_v=\Delta
v,\;\nabla\cdot
v=0,\mbox{ in }(0,T)\times\OO,\\ v=0\mbox{ on }(0,T)\times\po,\quad
v=v_\circ\mbox{ on
}\{0\}\times\OO.\ea\ee In system
\rf{NS}  $v$ is the kinetic field,
$\pi_v$ is the pressure field. We set 
 $b_t:=
\frac\partial{\partial t}b$  and
 $b\cdot\nabla d:=
b_k\frac{\partial d}{\partial x_k}$. 
In order to highlight the main ideas we assume:  $\OO\subseteq\R^3$ smooth bounded or exterior domain,  zero body
force and homogeneous boundary data.\par The symbol
$\mathscr C_0(\OO )$ stands for
the subset of $C_0^\infty(\OO )$
whose elements are divergence
free. We set
$J^2(\OO)\!:=$completion of
$\mathscr C_0(\OO)$ with respect to the $L^2$-norm,
and $J^{1,2}(\OO)\!:=$ completion of
$\mathscr C_0(\OO)$  with respect to  the  $W^{1,2}(\OO)$-norm.
 We set
$(u,g)_D:=\intl{\!\!D
} u\cdot gdx$, and in the case of $D\equiv \OO$ we drop the subscript $D$. \par Following  Prodi , \cite{P}, we set 
 \begin{defi}\label{WS}{\sl Assuming $v_\circ\in J^2(\OO)$, a field   $v:(0,\infty)\times\OO \to\R^3$   is said a weak solution to problem {\rm\rf{NS}} if \begin{itemize}
 \item[\rm i)] for all $T>0$,
 $v\in L^\infty(0,T;L^2(\OO))\cap L^2(0,T; J^{1,2}
(\OO ))$,  
\item[\rm ii)] $\displ\lim_{t\to0}\dm
v(t)-v_\circ\dm_2=0$, \item[\rm iii)] for all $t,s\in(0,T)$  the field $v$
satisfies the  equation:
\newline $$\displ\intll
st\Big[(v,\varphi_\tau)-(\nabla
v,\nabla
\varphi)+(v\cdot\nabla\varphi,v) \Big]d\tau+(v(s),\varphi
(s))=(v(t),\varphi(t)),$$
\newline\centerline{for all $  \varphi\in \mathscr W(\OO_T)$}, \end{itemize}where   the test functions set is defined as $$\ba{c}\mathscr W(\OO_T):=\{ \varphi\in C(0,T;J^{1,2}(\OO )),\mbox{ with } \varphi_t\in L^2(0,T;L^2(\OO)), \VS\hskip1.5cm\mbox{and } \varphi=0 \mbox{ in neighborhood of }T\}(PUNTO).\ea$$}
\end{defi} The following existence result holds:\begin{tho}\label{EXT}{\sl For any $v_\circ\in J^2(\OO)$ there exists a weak solution to problem \rf{NS} such that
\be\label{EI}\dm v(t)\dm_2^2+2\intll 0t\dm \nabla v(\tau)\dm_2^2d\tau\leq \dm v_\circ\dm_2^2, \mbox{ for all }t> 0,\ee and $(v(t),\psi)\in C([0,T))$ for all $\psi \in J^2(\OO)$.}\end{tho} The above existence result is due to  Hopf in \cite{H}. Inequality \rf{EI} is called energy inequality in weak form. It is different from the  following one   \be\label{EIL}\dm v(t)\dm_2^2+2\intll st\dm \nabla v(\tau)\dm_2^2d\tau\leq \dm v(s)\dm_2^2, \mbox{ for all }t>s,\mbox{ a.e.  in } s>0\mbox{ and for }s=0,\ee called energy inequality in strong form, and from the following one, which is a localized form of \rf{EI} and \rf{EIL}, that we   state for the Cauchy problem\be\label{SEI}\ba{l}\displ\intl{\R^3}|v(t)|^2\phi(t)dx+2\intll st\intl
{\R^3}|\nabla v(\tau)|^2\phi\, dxd\tau\leq \intl{\R^3}|v(s)|^2\phi(s)dx\VS\hskip 2,5cm+\intll st\intl{\R^3}|v|^2(\phi_\tau+\Delta\phi)dxd\tau+\intll st\intl{\R^3}(|v|^2+2\pi_v)
v\cdot\nabla\phi dxd\tau,\ea\ee for all $t\geq s$, for $s=0$ and a.e. in $s> 0$, and for all 
nonnegative $\phi\in C_0^\infty(\R\times\R^3)$. Inequality \rf{EIL} is due to  
Leray in \cite{L}  for the Cauchy problem.  
Subsequently, in the case of an IBVP, we have    the energy inequality in strong form for solutions   in $\OO$ bounded  considering again a  Hopf weak solution but   constructed  by means of  the   Heywood's device \cite{Hy}. 
    For exterior domains and more in general for unbounded domains there are several contributions, see for example \cite{MS,FK}.  Inequality \rf{SEI} is due to Caffarelli, Khon and Nirenberg in \cite{CKN}. To date, in unbounded three dimensional domains the energy inequality in strong form and in localized form is not proved   for the solutions  furnished  by Hopf's technique, \cite{H}. Moreover, regardless of the kind of weak solution,  it is not known if  the energy inequality in strong form holds for $\Omega\subset \R^n$,  $n>4$ and $\Omega$ unbounded domain. However  the result claimed in Theorem\,\ref{EXT}, or its variant in \cite{L} and in \cite{CKN}, is the unique existence result at disposal for arbitrary data.
 \par  
 It is known that the regularity of a weak solution to   problem \rf{NS} is an open question (see e.g. \cite{L,Ld,CKN}). In the interval     between the two essays  \cite{L} and \cite{CKN},     Prodi and  Serrin, in the papers \cite{P} and \cite{Sr0,Sr1}, respectively,   introduce   the idea of  searching for sufficient conditions in order to obtain the energy equality, the uniqueness and the regularity of a weak solution. This approach translates into   extra assumptions     that are apt to obtain the well posedeness of the problem. A well known and classical result concerns the  regularity and the uniqueness: \par\centerline{\sl if $v$ is a weak solution and $v\in L^\rho(0,T;L^\sigma(\OO))$, then $v$ is    smooth and unique,}\par\noindent provided that   $\frac 3\sigma+\frac2\rho=1$ with $\sigma\in(3,\infty]$, $\rho\in[2,\infty)$\,. A  proof of this result is given by Giga in \cite{Gy} both in the cases of the Cauchy problem and IBVP with $\OO$ bounded. In arbitrary domains the problem is considered by Galdi and Maremonti in \cite{GM}. The limit case $L^\infty(0,T;L^3(\OO))$ is not   considered   and is studied by Escaurazia, Seregin and \v{S}ver\'{a}k  in \cite{ESS}.  In \cite{KK}, Kim and Kozono establish interior regularity considering    Lorentz spaces in place  of   Lebesgue spaces. In this connection a further contribution is given  by Bosia, Pata and Robinson in \cite{BPR}.   There exists a wide literature on  sufficient conditions for the regularity of a weak solution.  A possible key tool  to obtain these results   is the one based on the mild solutions to the integral equation associated to problem \rf{NS}. They are established on the wake of the ones   due to Kato in \cite{KT} and to Giga in \cite{Gy}.  In this connection see the paper \cite{F}, where Farwig provides  an interesting review of the state of the art on the  problematic and on the techniques. Different assumptions, closely connected with the one by Prodi and by Serrin,   are considered in a series of papers.    An interesting   update on the topic is given    in the recent paper  by Tran and Yu \cite{TY}.
 \par In \cite{CKN}, a new highlight  on the Prodi-Serrin ideas is given by Caffarelli, Khon and Nirenberg. They detect that a solution $v$ and the  condition $v\in L^\rho(s,t;L^\sigma(D))$ ($\frac3\sigma+\frac2\rho=1$) suitably satisfy a  scaling invariance property.  That is if $v$ is a solution    then for all $\lambda>0$ we get that $v_\lambda(t,x):=\lambda v(\lambda^2t,\lambda x)$ is still a solution, and $\dm v_\lambda\dm_{L^\rho(s,t;L^\sigma( D))}\equiv\dm v\dm_{L^\rho(s,t;L^\sigma(D))}$ independently of the domain $(s,t)\times D$.   It is just the case to recall that the kind  of scaling invariant norm is connected with dimensional balance of   equation \rf{NS}$_1$. In \cite{CKN} by means of Proposition\,1 and Proposition\,2 the authors realize  local   (that is on a space-time parabolic neighborhood of  points)  estimates in  norms which are scaling invariant. Then they obtain  new sufficient conditions for the regularity  proving partial regularity for a suitable weak solution  with an initial data in $J^2(\OO)$. In this connection,   starting from an idea already contained in \cite{CKN}, in the recent paper \cite{CM} Crispo and Maremonti detect that the bound of the metrics employed in \cite{CKN} for a weak solution is ensured by some weighted norms, whose advantage is   that they hold for all $t>0$ provided that the initial data  satisfies suitable assumptions.\par  Hence the possibility of realizing the regularity seems connected with the existence of   scaling invariant norms of a weak solution $v$ on  $(0,T)\times\OO$.  \par The above considerations on the scaling invariant metrics  have generated, explicitly   or tacitly, a way of thinking for which in order to obtain regularity for a weak solution also   the initial data $v_\circ$ has to belong to some function space  which is scaling invariant.     As a matter of course this leads to consider initial data in more regular spaces than $J^2(\OO)$.    Conversely, it is natural to inquire about the compatibility between an initial data {\it a priori} in $J^2(\OO)$ and the regularity of solutions for $t>0$. In other words one questions if it is well posed  the following \begin{defi}\label{RS}{\sl  We say  that a weak solution $u$   is a regular solution to problem \rf{NS} if for all $\vep>0$ and $T>\vep$ we have  $u\in L^\infty((\vep,T)\times \OO)$\,.}\end{defi} 
An analogous question can be posed on the extra assumption related to a weak solution $v$ of Theorem\,\ref{EXT} in order to obtain the energy equality:
\be\label{PEEO}\dm v(t)\dm_2^2+2\intll st\dm \n v(\tau)\dm_2^2d\tau =\dm v(s)\dm_2^2,\;\mbox{ for all }t>s\geq0.\ee In \cite{P},   under the extra assumption   $v\in L^4(0,T;L^4(\OO))$ Prodi proves the energy equality. More recently, in \cite{CFS} and \cite{CCFS} the assumptions are different.  In \cite{CFS} Cheskidov, Friedlander and Shvydkoy assume $\intll0T\dm A^\frac5{12}u(t)\dm_2^3dt$, where $A$ is the Stokes operator. In \cite{CCFS} Cheskidov, Costantin, Friedlander and Shvydkoy consider the assumption $v\in L^3(0,T;B^\frac13_{3,\infty}(\R^3))$, where $B^\frac13_{3,\infty}$ is the Besov space. Finally\,\footnote{\; For more general domains see also \cite{FT}.} in \cite{F} Farwig assumes $\intll0T\dm A^\frac14v(t)\dm_\frac{18}7^3dt$. These assumptions have the same scaling, in particular we get $v\in L^3(0,T;L^\frac92(\OO))$ which furnishes $\frac3{\frac92}+\frac23=\frac43>\frac 34+\frac12$, the last one being the  Prodi condition\,\footnote{\,  Concerning the energy equality, also in connection with the one in local form \rf{SEI}, an interesting analysis can be found in \cite{JM}.}. In any case these conditions are not compatible on $(0,T)$ with an initial data in $J^2(\OO)$. \par
The main goal of this note is to prove the following results related to the regularity and to the energy equality, and their compatibility with the initial data in $J^2(\OO)$.
  \begin{tho}\label{CT}{\sl Assume that the weak solution $v$ of Theorem\,{\rm\ref{EXT}} satisfies the  condition:\be\label{EX}\mbox{for all }\vep>0\;v\in L^\rho(\vep,T;L^\sigma(\OO))\mbox{ with $\frac3\sigma+\frac2\rho=1$,}\;\rho<\infty,\ee then $v$ is a regular solution, and 
\be\label{LLCR-II}(t-\vep) \dm v_t\dm_2\leq  B(\dm v_\circ\dm_2,\vep,t) ,\mbox{ for all }t>\vep\,,\ee where   $ B(\dm v_\circ\dm_2,\vep,t):=c\dm v_\circ\dm_2  \exp[t-\vep+c\intll\vep t\dm v(\tau)\dm_\sigma^\rho d\tau] $, and $c$ is a constant independent of $v$ and $\vep$.}\end{tho}
For the energy equality we do not consider the   conditions furnished in \cite{CFS,CCFS,F}. We limit ourselves to prove
\begin{tho}\label{CTO}{\sl Assume that the weak solution $v$ of Theorem\,{\rm\ref{EXT}} satisfies the  condition:\be\label{CTO-I}\mbox{for all }\vep>0,\;v\in L^4(\vep,T;L^4(\OO)), \ee then for  $v$  the energy equality \rf{PEEO} holds.
  In particular we get that $v\in C([0,T);J^2(\OO))$.}\end{tho}
Assumptions  \rf{EX} and \rf{CTO-I} are a weak form of Prodi-Serrin conditions that yield  the analogous results.  \par By interpolation of the spaces $L^\infty(0,T;L^2(\OO))$ and $L^\infty((\vep,T)\times\OO)$ all  regular solutions belong to $L^\rho(\vep,T;L^\sigma(\OO))$ with $\frac3\sigma+\frac2\rho=1$, therefore the set of regular solutions is characterized by means of the extra condition \rf{EX}.\par Via propetry \rf{LLCR-II} a regular solution is a classical solution for $t>0$ (see e.g. \cite{Sr0} and also \cite{O}). \par We remark that, setting $\mathbb L^\sigma$:=completion of $\mathscr C_0$ in the Lorentz space $ L(\sigma,\infty)$, we can consider the assumption $v\in L^\rho(\vep,T;\mathbb L^\sigma(\OO))$ which is a weaker spatial assumption, close to the one employed in \cite{KK}.  It is important to point out that, by interpolation, for $\sigma\in[3,6]$   the assumption $v\in \mathbb L^\sigma$ is automatically satisfied. \par We stress that for a weak solution of Theorem\,\ref{EXT} condition \rf{EX}   furnishes the regularity in the sense of Definition\,\ref{RS},  but we are not able to prove that assumption \rf{EX} also implies uniqueness.  This makes the difference with the Prodi-Serrin condition  which ensures  both the properties.  
\par {\it Mutatis mutandis} the notion of regular solution can be also given for solutions $u$ to the Stokes problem (that is \rf{NS}   dropping the convective term). In this case, if $u_\circ\in L^p(\OO)$, then it is possible to give a behavior in $t=0$  of the $L^q$-norm of the solutions, $q\geq p$. Actually, we get $\displ \lim_{t\to0}t^{\frac32(\frac1p-\frac1q)}\dm u(t)\dm_q=0$ and $\dm u(t)\dm_q\leq c\dm u_\circ\dm_pt^{-\frac32(\frac1p-\frac1q)}$, with $c$ independent of $u_\circ$. The case of data in $L(p,\infty)(\OO)$ is different. There the above $L^p-L^p$ properties are admissible  for initial data $u_\circ \in \mathbb L(p,\infty)$ (see e.g. \cite{Mvi}).
 In the two dimensional case  the definition of regular solution is well posed (see Ladyzhenskaya \cite{Ld}). In this case we also get  behavior in $t=0$ of the solution for $q\geq p=2$. This is consequence of the estimate $\dm \nabla v(t)\dm_2\leq \dm v_\circ\dm_2\exp(c\dm v_\circ\dm_2)t^{-\frac12} $, the energy equality and the Gagliardo-Nirenberg inequality. In \cite{KY}, Kozono-Yamhazaki furnish a generalization of Ladyzhenskaya's result. They assume  $v_\circ\in L(2,\infty)$ and divergence free, but in this case no result about the behavior of the solution $v$ in a neighborhood of $t=0$ is known. In particular the estimate for $q=p=2$ does not hold. We stress that in all the  listed cases   of  regular solutions $u$ the uniqueness holds.   \par The plain of the note is the following. In sect. 2 we prove Theorem\,\ref{CTO}. This is done by means of the propreties of solutions to a suitable linearized Navier-Stokes problem. In particular we deduce the energy equality. Out of respect for G. Prodi, our proof is developed following in the first step  the argument lines given in \cite{P}. Sect.\,3 is devoted to some auxiliary results. Finally, in sect.\,4 we give the proof of Theorem\,\ref{CT}. In doing this we partially follow the argument lines by Galdi in \cite{G}, that   follows in turn the ones by Galdi and Maremonti in \cite{GM}.   
\section{An improvement of the Prodi  result: the energy equality.}
Fundamental for our aims is the study of the linearized Navier-Stokes problem:
\be\label{LNS}\ba{l}w_t+a_r\cdot
\nabla w+\nabla\pi_w=\Delta
w,\;\nabla\cdot
w=0, \mbox{ in }(0,T)\times\OO,\\ w=0\mbox{ on }(0,T)\times\po,\quad
w=w_0\mbox{ on
}\{0\}\times\OO,\ea\ee  where for the coefficient $a_r$ the subscript  $r$  ranges between $1$ or $2$. Respectively, the coefficients  enjoy the following integrability properties: \be\label{AA}\ba{ll}&a_1 \in L^2(0,T;J^{1,2}(\OO))\cap L^4(\vep,T;L^4(\OO)),\VSE a_2 \in L^2(0,T;J^{1,2}(\OO))\cap L^\rho(\vep,T;L^\sigma(\OO)),\mbox{ $\frac3\sigma+\frac2\rho=1$}.\ea\ee Roughly speaking,  assuming the coefficients only  in $L^2(0,T;J^{1,2}(\OO))$   the weak solutions to the linearezed problem \rf{LNS} reveal the same difficulties in order to prove either an energy equality relation and uniqueness of the solutions. We are interested in both the questions, because we reduce the proof of Theorem\,\ref{CT} to the study of a suitable linearezed problem. For these aims it is also crucial to consider  the mollified system
\be\label{LNSM}\ba{l}w_t+J_n[a]\cdot
\nabla w+\nabla\pi_w=\Delta
w,\;\nabla\cdot
w=0,\mbox{ in }(0,T)\times\OO,\\ w=0\mbox{ on }(0,T)\times\po,\quad
w=w_0\mbox{ on
}\{0\}\times\OO,\ea\ee where, for all $n\in\N$,  $J_n[\cdot]$ is the time-space  Friderichs mollifier, and $a\in L^2(0,T;J^{1,2}(\OO))$, and it is extended to $0$ for $t<0$. \par 
The following result holds (cf. Solonnikov \cite{Sl}):
\begin{tho}\label{LNSS}{For all $\phi_0\in J^{1,2}(\OO)$ there exists a unique solution $(\phi,\pi_\phi)$ to problem \rf{LNSM} such that $\phi\in C(0,T;J^{1,2}(\OO))\cap L^2(0,T;W^{2,2}(\OO))$ and $\phi_t,\,\n\pi_\phi\in L^2(0,T;L^2(\OO))$. Moreover, uniformly with respect to $n\in\N$, the energy equality holds:
\be\label{EEW}\dm \phi(t)\dm_2^2+\intll st\dm\n \phi(\tau)\dm_2^2d\tau=\dm \phi(s)\dm_2^2,\mbox{ for all }t>s\geq 0.\ee}\end{tho} In order to work with the weak solutions to problem \rf{NS}, it is better to study a weak formulation of problem \rf{LNS}. 
\begin{defi}\label{LNSWS}{\sl 
Assuming $w_\circ\in J^2(\OO)$, a field   $w:(0,\infty)\times\OO \to\R^3$   is said a weak solution to problem {\rm\rf{LNS}} if \begin{itemize}
 \item[\rm i)] for all $T>0$,
 $w\in L^\infty(0,T;L^2(\OO))\cap L^2(0,T; J^{1,2}
(\OO ))$,  
\item[\rm ii)] $\displ\lim_{t\to0}\dm
w(t)-w_\circ\dm_2=0$, \item[\rm iii)] for all $t,s\in(0,T)$  the field $w$
satisfies the  equation:
\newline $$\displ\intll
st\Big[(w,\varphi_\tau)-(\nabla
w,\nabla
\varphi)+(a\cdot\nabla\varphi,w) \Big]d\tau+(w(s),\varphi
(s))=(w(t),\varphi(t)),$$ for all $\varphi\in \mathscr W(\OO_T)$.\end{itemize}}
\end{defi}
  In \cite{P}, for all weak solutions to problem \rf{NS} Prodi proves the energy equality, that is
$$\dm v(t)\dm_2^2+2\intll0t\dm \n v(\tau)\dm_2^2d\tau=\dm v_\circ\dm_2^2,\mbox{ a.e. in }t>0,$$
provided that they enjoy  the extra condition\,\footnote{Actually, as  was recognized and remarked subsequently in the literature, by means of the   extra assumption $v\in L^4(0,T;L^4(\OO))$, the result of the energy equality has  $n$-dimensional validity, $n\geq2$.} $v\in L^4(0,T;L^4(\OO))$. Actually,   Prodi's result contains inside a uniqueness theorem for weak solutions to problem \rf{LNS}.   The following lemma and related corollary give an improvement of the quoted results because the assumption is relaxed as follows:   $L^4(\vep,T;L^4(\OO))$ for all $\vep>0$.
\begin{lemma}\label{EEI}{\sl Assume that $w$ is a weak solution to problem \rf{LNS} with coefficient $a\equiv a_1$. Assume that, for all $\vep>0$, $w\in L^4(\vep,T;L^4(\OO))$, $(w(t),\psi)\in C((0,T))$ for all $\psi\in J^2(\OO)$, and the energy inequality (in weak  form) holds: 
\be\label{EIS}\dm w(t)\dm_2^2+2\intll 0t\dm \n w(\tau)\dm_2^2d\tau\leq \dm w_0\dm_2^2,\mbox{ for all }t>0 .\ee     Then the energy equality holds for $w$: \be\label{EE}\dm w(t)\dm_2^2+2\intll st\dm \n w(\tau)\dm_2^2d\tau=\dm w(s)\dm_2^2,\mbox{ for all }t>s\geq0 .\ee  }\end{lemma}
\Pr The proof of property \rf{EE} is achieved in two steps. In the first step we employ  the Prodi technique. We denote by $N$ the set of those $t$ such that $\dm \n w(t)\dm_2<\infty$ and consider $s\in(0,t)$. We define   \be\label{DWS}w^*(\tau,x):=w(\tau,x)- w(t,x)\frac{\tau-s}{t-s}\mbox{ for }\tau\in[s,t]\mbox{ and }w^*(\tau,x)=0\mbox{ for }\tau\notin[s,t]\,.\ee Moreover we set $$w_\eta(\tau,x)=h(\tau,t)w(t)\frac{\tau-s}{t-s} + J_\eta[J_\eta[w^*]](\tau)\,,$$ where the function $h(\tau,t)$ is a nonnegative smooth cutoff function such that $h(\tau,t)=1$ for $\tau\leq t$ and $h(\tau,t)=0$ for $\tau\geq2t$, and   $J_\eta$ is a mollifier. It is easy to check that $w_\eta\in \mathscr W(\OO_T)$. So that we use $w_\eta$ as test function in  iii) of Definition\,\ref{LNSWS}:\be\label{IEEE}\intll
st\Big[(w,{w_\eta}_\tau)-(\nabla
w,\nabla
w_\eta)+(a_1\cdot\nabla w_\eta,w)\Big]d\tau+(w(s),w_\eta
(s))=(w(t),w_\eta(t))\,.\ee
We evaluate the terms:
$$(w(t),w_\eta(t))-\intll
st(w,{w_\eta}_\tau) d\tau =:I_1+I_2\,.$$ We get
$$I_1=\dm w(t)\dm_2^2+(w(t),J_\eta[J_\eta[w^*]](t)).$$ By virtue of the definition \rf{DWS} of  $w^*$, making use of an integration by parts, we get
$$\ba{ll}I_2\hskip-0.2cm&\displ=-\intll st(w(\tau),w(t)(t-s)^{-1}+\mbox{$\frac{\partial}{\partial \tau}$}J_\eta[J_\eta[w^*]](\tau))d\tau\VSE=-\intll st(w^*(\tau)+  w(t)\frac{\tau-s}{t-s},w(t)(t-s)^{-1}+\mbox{$\frac{\partial}{\partial \tau}$}J_\eta[J_\eta[w^*]](\tau))d\tau\VSE=-\mbox{$\frac12$}\dm w(t)\dm_2^2-\intll st(w^*(\tau),w(t)(t-s)^{-1}+\mbox{$\frac{\partial}{\partial \tau}$}J_\eta[J_\eta[w^*]](\tau))d\tau\VSE\hskip3cm-(w(t),J_\eta[J_\eta[w^*](t))+(t-s)^{-1}\!\!\intll st(w(t), J_\eta[J_\eta[w^*]](\tau))d\tau\,.\ea$$ Summing we get
$$\ba{ll}I_1\!+I_2\hskip-0.3cm&=\frac12\dm w(t)\dm_2^2-\displ\!\!\intll st\! (w^*(\tau),\mbox{$\frac{\partial}{\partial \tau}$}J_\eta[J_\eta[w^*]](\tau))d\tau\!+\!\mbox{$\frac1{t-s}$}\!\!\intll st \!(w(t), J_\eta[J_\eta[w^*]](\tau)\!-w^*(\tau))d\tau\VS& =\frac12\dm w(t)\dm_2^2-H_1(\eta)+H_2(\eta)\,.\ea$$Recalling the definition of $w^*$, we obtain the identity
$$\ba{ll}H_1(\eta)\hskip-0.2cm&\displ=\intll{-\infty}t(w^*(\tau), \intll{-\infty}{\infty}\mbox{$\frac{\partial}{\partial \tau}$}J_\eta(\tau-h)\intll{-\infty}\infty J_\eta(h-s)w^*(s)dsdh)d\tau\VSE=-\intll{-\infty}\infty(\mbox{$\frac{\partial}{\partial h}$}\intll{-\infty}\infty J_\eta(h-\tau)w^*(\tau)d\tau,\intll{-\infty}{\infty}J_\eta(h-s)w^*(s)ds)dh\VSE=-\mbox{$\frac12$}\intll{-\infty}\infty\mbox{$\frac{d}{dh}$}\dm J_\eta[w^*](h)\dm_2^2dh=0, \mbox{ for all }\eta>0\,.\ea$$ Now we evaluate the limit in $\eta\to0$ of $H_2(\eta)$. Since $J_\eta[J_\eta[w^*]]\to w^*$ in $L^2(0,T;L^2(\OO))$, we have
 $$\lim_{\eta\to0}|H_2(\eta)|\leq\lim_{\eta\to0} \mbox{$\frac1{t-s}$}\intll st\dm w(t)\dm_2\dm J_\eta[J_\eta[w^*]](\tau)-w^*(\tau)\dm_2d\tau=0.$$ 
For the term $(w(s),w_\eta
(s))$ we get
$$\ba{ll}(w(s),w_\eta
(s))\hskip-0.2cm&\displ=(w(s),\intll{-\infty}\infty J_\eta(s-h)\intll st J_\eta(h-\xi)w^*(\xi)d\xi dh)\VSE=\mbox{$\frac12$}\dm w(s)\dm_2^2+(w(s),\intll{-\infty}\infty J_\eta(s-h)\intll{s}t J_\eta(h-\xi)(w^*(\xi)-w(s))d\xi dh).\ea$$ Since $(w(s),w^*(\xi))\in C(s,t)$, and $\displ\lim_{\xi\to s }(w(s),w^*(\xi)-w(s))=0$,   we get
$$\lim_{\eta\to0}(w(s),w_\eta(s))=\mbox{$\frac12$}\dm w(s)\dm_2^2. $$ Hence the limit for $\eta\to0$   gives
$$\lim_{\eta\to0}\Big[(w(t),w_\eta(t))-\!\!\intll
st(w,{w_\eta}_\tau) d\tau\! -(w(s),J_\eta[J_\eta[w(s)]])\Big]\!= \mbox{$\frac 12\dm w(t)\dm_2^2-\frac12\dm w(s)\dm_2^2$}\,.$$ We consider the weak limit of $J_\eta[J_\eta[w^*]]$ in $L^2(0,T;J^{1,2}(\OO))$. Hence recalling definition \rf{DWS} of $w^*$, we get
$$\ba{ll}\displ\lim_{\eta\to0}\intll st(\nabla w,\nabla w_\eta)d\tau\hskip-0.2cm&\displ=\lim_{\eta\to0}\intll st(\nabla w,\n w(t)\mbox{$\frac{\tau-s}{t-s}$}+\nabla J_\eta[J_\eta[w^*]])d\tau\VSE=\lim_{\eta\to0}\intll st(\nabla w,\nabla J_\eta[J_\eta[w]])d\tau=\intll st\dm \nabla w(\tau)\dm_2^2d\tau\,.\ea$$ Finally, we evaluate the limit of the nonlinear part. In this limit we employ the assumption of $a_1\in L^4(\vep,T;L^4(\OO))$. Actually, employing the fact that  $(a_1\cdot\nabla w,w)=0$ almost everywhere in $t>0$, recalling formula\rf{DWS} and $s-2\eta>\frac s2$, we get
$$\displ \lim_{\eta\to0}\intll st\!(a_1\cdot\nabla w,w_\eta)d\tau\displ=\lim_{\eta\to0}\intll st\!(a_1\cdot\nabla w,w_\eta-w)d\tau=\lim_{\eta\to0}\intll st\!(a_1\cdot\n w, J_\eta[J_\eta[\chi_{[s,t]} w]]-w)d\tau,$$ where $\chi$ is the charateristic function of the interval $(s,t)$. Applying H\"older's inequality to last term, we obtain
$$\ba{ll}\displ|\lim_{\eta\to0}\intll st(a_1\cdot\n w,w_\eta)d\tau|\hskip-0.2cm&\displ\leq\lim_{\eta\to0}\intll st\dm a_1\dm_4\dm\nabla w\dm_2\dm J_\eta[J_\eta[\chi w]]-w\dm_4d\tau\VSE\leq\lim_{\eta\to0} \Big[\intll st\!\dm a_1\dm_4^4d\tau\Big]^\frac14\Big[\intll st\! \dm\nabla w\dm_2^2d\tau\Big]^\frac14 \Big[\intll st \!\dm J_\eta[J_\eta[\chi w]]-w\dm_4^4d\tau\Big]^\frac14\!.\ea$$ Hence we get
$$\lim_{\eta\to0}|\intll st(a_1\cdot\nabla w,w_\eta)d\tau| =0.$$ 
Considering  each limit for $\eta\to0$ for the corresponding term of \rf{IEEE}, we deduce  
\be\label{EISI}\dm w(t)\dm_2^2+2\intll st\dm \n w(\tau)\dm_2^2d\tau=\dm w(s)\dm_2^2,\mbox{ for all }t\in N,\mbox{ and }  s>0 .\ee 
Since \rf{EISI} holds for all $s>0$, by virtue of ii) of Definition\,\ref{LNSWS} and the absolute continuity of the integral function, from the above equality we deduce
\be\label{SEE}\dm w(t)\dm_2^2+2\intll 0t\dm \n w(\tau)\dm_2^2d\tau=\dm w_0\dm_2^2,\mbox{ for all }t\in N.\ee Now we prove the property for all $t>0$. To this end we prove that for all $t>0$   $\displ\lim_{\xi\to t^+}\dm u(\xi)-u(t)\dm_2=0$ holds. We consider problem \rf{LNSM} with coefficient $J_n[\widehat a]$, with   $\widehat a:=a_1(\xi-h,x)$ for $h\in[t,\xi]$ otherwise $\widehat a=0$ for $h\notin[t,\xi]$. We have  $\widehat a\in L^2(0,\xi;J^{1,2}(\OO)) $. By virtue of Theorem\,\ref{LNSS}, for all $n\in\N$, we obtain the solution $(\phi_n,\pi_{\phi_n})$. We set $\widehat\phi_n(\tau,x):=\phi_n(\xi-\tau,x)$, for all $\tau\in[t,\xi].$ Since $\widehat\phi_n$ is solution backward in time, substituting $\widehat\phi_n$ in i) of Definition\,\ref{LNSWS}, and integrating by parts on $(t,\xi)\times\OO$, we get $$(w(\xi),\phi_0)=(w(t),\phi_n(\xi-t))+\intll t\xi((J_n[a_1](\tau) -a_1(\tau))\cdot\n\widehat\phi_n,w)d\tau,\mbox{ for all }n\in\N.$$ On the other hand, for all $\xi>t>0$ and uniform in $n\in\N$ estimate \rf{EEW} ensures
$$\dm \phi_n(\xi-t)\dm_2\leq\dm \phi_0\dm_2\mbox{ and }\intll t\xi\dm \n \phi_n(t-\tau)\dm_2^2d\tau\leq \dm \phi_0\dm_2^2\,.$$ Therefore, applying H\"older's inequality, we deduce $$\ba{ll}|(w(\xi),\phi_0)|\hskip-0.3cm&\displ\leq \dm w(t)\dm_2 \dm \phi_0\dm_2+\!\Big[\!\intll t\xi\!\dm J_n[a_1](\tau) -a_1\dm_4^4\Big]^\frac14\Big[\!\intll0{\xi-t}\!\!\dm \n\phi(h)\dm_2^2 dh\Big]^\frac14\Big[\intll t\xi\!\dm w(\tau)\dm_4^4d\tau\Big]^\frac14\\&\displ\leq \dm\phi_0\dm_2\Big[ \dm w(t)\dm_2 + \Big[\!\intll t\xi\dm J_n[a_1](\tau) -a_1\dm_4^4\Big]^\frac14 \Big[\intll t\xi\dm w(\tau)\dm_4^4d\tau\Big]^\frac14\Big]\,\mbox{ for all }n\in\N.\ea$$ Since $\phi_0$ is arbitrary in $J^{1,2}(\OO)$, in the limit for $n\to\infty$ 
we deduce $$\dm w(\xi)\dm_2\leq\dm w(t)\dm_2\mbox{ for all }\xi>t>0.$$ This last property and the assumption of $(w(t),\psi)$ continuous function of $t$, ensure that $w$ is continuous in $t$ on the right  in $L^2$-norm, for all $t\geq 0$. Let $t\in (0,T)-N$. Then for all sequence $\{t_n\}$ which converges to $t$ from the right hand side  we have the limit property $\displ\lim_{t_n\to t}\dm w(t_n)\dm_2=\dm w(t)\dm_2$. Therefore from \rf{SEE} written for the instant of the sequence $\{t_n\}$ we deduce \rf{EE} for all $t>0$ and $s=0$. After that easily follows \rf{EE} for all $t>s\geq0$.\chiu The following result immediately holds:\begin{coro}\label{C-I}{\sl In the hypotheses of Lemma\,{\rm\ref{EEI}}, we get: $$\mbox{for all }T>0,\; w\in C([0,T;J^2(\OO)),$$ and  if,  for some $s\geq0$, $\dm w(s)\dm_2=0$ then $w$ is identically null for all $t>s$.}\end{coro}The following result proves Theorem\,\ref{CTO}:\begin{coro}\label{PEE}{Let $v$ be a weak solution to problem \rf{NS}. Assume that for all $\vep>0$ we have $v\in L^4(\vep,T;L^4(\OO))$, then for $v$ the energy equality holds:
\be\label{PEE-I}\dm v(t)\dm_2^2+2\intll st\dm \n v(\tau)\dm_2^2=\dm v(s)\dm_2^2,\;t>s\geq0\,.\ee}\end{coro}\Pr It is enough to apply Lemma\,\ref{EEI} considering $v$ as a weak solution to problem \rf{LNS} with $a_1:=v$.\chiu
\section{Some   auxiliary results} We recall some well known estimates.
\begin{lemma}\label{GN}{\sl Assume that $D^2v\in L^p(\OO)$ and $v\in W_0^{1,r}(\OO)$. Then there exists a constant $c$ independent of $v$ such that
\be\label{GNI}\ba{lll}\dm \n v\dm_q&\hskip-0.33cm\leq c\dm D^2v\dm_p^\lambda\dm \n v\dm_r^{1-\lambda}, &\mbox{ provided that }\mbox{$\frac 1q=\lambda(\frac 1p-\frac1n)+(1-\lambda)\frac1r\,,$}\VS \hskip0.33cm\dm v\dm_\infty&\hskip-0.35cm\leq c\dm D^2v\dm_p^{\ov\lambda} \dm v\dm_r^{1-\ov{\lambda}},&\mbox{ provided that }\mbox{$0=\ov\lambda(\frac1p-\frac2n)+(1-\ov\lambda)\frac1r$}\,.\ea\ee}\end{lemma} 
\Pr For $\OO$ exterior domain estimates \rf{GNI}   are particular cases of the one proved in \cite{CM}. For $\OO$ bounded domain  or $\OO\equiv\R^n$, inequality \rf{GNI} is a particular case of the well known Gagliardo-Nirenberg inequality.  \chiu
\begin{lemma}\label{TR}{\sl Assume that $v\in W^{2,2}(\OO)\cap J^{1,2}(\OO)$, then 
\be\label{SD}\dm D^2v\dm_2\leq c(\dm P\Delta w\dm_2+\dm v\dm_{L^2(D)}),\ee where $D\subseteq\OO$ is a bounded domain  with $\partial(\OO-D)\cap \po=\emptyset$, and the constant $c$ is independent of $v$.}\end{lemma}
\Pr See \cite{H}. 
\begin{lemma}\label{UTR}{\sl Let $v$ be as in Lemma\,{\rm\ref{GN}}. Assume that $u\in L^2(\OO)$ and, for $q\in [2,6)$, $b\in L^{\frac{2q}{q-2}}(\OO)$, then we have
\be\label{TRI}|(b\cdot\n v,u)|\leq \mbox{$\frac 14$}\dm P\Delta v\dm_2^2+\mbox{$\frac14$}\dm u\dm_2^2+ c(\dm b\dm_{\frac{2q}{q-2}}^2+\dm b\dm_{\frac{2q}{q-2}}^\frac {4q}{6-q})\dm \n v \dm_2^2 \,.\ee   Moreover, if $w\in J^{1,2}(\OO)$ and $q>3$, then we have
\be\label{TRII}|(w\cdot\n w,b)|\leq c\dm w\dm_2^2\dm b\dm_{q}^\frac{2q}{q-3}+\mbox{$\frac12$}\dm \n w\dm_2^2\,.\ee In inequalities \rf{TRI} and \rf{TRII} the constant $c$ is independent of $v,u,b$ and $w$.}\end{lemma} \Pr In the case of \rf{TRI}, applying H\"older's inequality, we get
$$|(b\cdot\n v, u)|\leq \dm b
\dm_{\frac{2q}{q-2}}\dm \n v\dm_q\dm u\dm_2\,.$$   Applying  estimate \rf{GNI} with $p=r=2$, and subsequently \rf{SD} we obtain 
$$\ba{ll}|(b\cdot\n v, P\Delta v)|\hskip-0.2cm&\leq \dm b
\dm_{\frac{2q}{q-2}}\dm \n v\dm_q\dm u\dm_2\leq c\dm b
\dm_{\frac{2q}{q-2}}\dm \n v\dm_2^{1-\lambda}\dm D^2v\dm_2^\lambda\dm u\dm_2\VSE\leq c\Big[\dm b
\dm_{\frac{2q}{q-2}}\dm \n v\dm_2^{1-\lambda}\dm P\Delta v\dm_2^{\lambda}\dm u\dm_2+\dm b
\dm_{\frac{2q}{q-2}}\dm \n v\dm_2^{1-\lambda}\dm u\dm_2\dm  v\dm_{L^2(D)}^{\lambda}\Big]\,,\ea$$ where the exponent $\lambda=\frac{3(q-2)}{2q}$. Finally, since $D$ is bounded and on $\po\cap\partial D$ we have $v=0$, by means of Poincar\é inequality, via the Cauchy inequality, we arrive at \rf{TRI}. Finally, for estimate \rf{TRII}, applying H\"older's inequality, and the Gagliardo-Nirenberg inequality    we easily obtain
$$|(w\cdot\n w,b)|\leq c\dm w\dm_{\frac{2q}{q-2}}\dm \n w\dm_2\dm b\dm_q\leq  c\dm w\dm_2^{\frac{q-3}q}\dm \n w\dm_2^{\frac {q+3}q}\dm b\dm_q.$$Hence  the Cauchy inequality leads to \rf{TRII}. \chiu  The following theorem concerns the existence of the so called strong regular solutions in the special case of the $L^2$-theory.
\begin{tho}\label{TR}{\sl Let $v_\circ\in J^{1,2}(\OO)$. The there exists a unique solution   to {\rm\rf{NS}} such that
\be\label{TR-I}v\in C(0,T;J^{1,2}(\OO))\cap L^2(0,T;W^{2,2}(\OO)),\mbox{ with }v_t,\n\pi\in L^2(0,T;L^2(\OO)),\ee and, for all $k\in\N,\,\ell=0,1,2,\,\vep>0$, \be\label{TR-II}\n^\ell D^k_tv\in C(\vep,T;L^2(\OO)).\ee}\end{tho}\Pr This result is a particular case of the one proved in Theorem\,3 in \cite{Hy}, see also Chap.\,V in  \cite{ES}.\chiu
\section{Proof of Theorem\,\ref{CT}} The proof of Theorem\,\ref{CT} is achieved by means of two lemmas.
\begin{lemma}\label{LCR}{\sl In the hypotheses of Theorem\,{\rm\ref{CT}} for all $\vep>0$ we get
\be\label{LCR-I}\ba{c}v\in C(0,T;J^2(\OO))\cap C(\vep,T;J^{1,2}(\OO))\cap L^2(\vep,T;W^{2,2}(\OO))\VS v_t,\n \pi_v\in L^2(\vep,T;L^2(\OO)),\ea\ee for arbitrary $T>\vep>0$. Moreover, we get
\be\label{LCR-II}\dm v(t)\dm_2^2+2\intll st\dm \n v(\tau)\dm_2^2d\tau=\dm v(s)\dm_2^2,\mbox{ for all }t>s\geq0 , \ee 
and \be\label{EGP}\ba{l}\displ \displ(t-\vep)\dm \n v(t)\dm_2^2+\mbox{$\frac12$}\intll\vep t(\tau-\vep)(\dm P\Delta v(\tau)\dm_2^2+\dm v_\tau(\tau)\dm_2^2)d\tau \VS \hskip 4cm  \leq \dm v_\circ\dm_2^2\exp[t-\vep+c\intll\vep t\dm v(\tau)\dm_\sigma^\rho d\tau]=:A^2(\dm v_\circ\dm_2,\vep,t) \,,\ea\ee with $c$   independent of $\vep$  and $t$.}\end{lemma}
\Pr We consider problem \rf{LNSM} with $a:=v$ and $w_\circ:=v_\circ^n$, where $\{v_\circ^n\}\subset\mathscr C_0(\OO)$ is a sequence which converges to $v_0$ in $J^2(\OO)$. By virtue of Theorem\,\ref{LNSS} we obtain a sequence $\{w^n\}$ of solutions to problem \rf{LNSM}. Now our task is to prove the existence of a limit $w$ which is a regular solution to problem \rf{LNS} with $a_2\equiv v$. We base the existence of the limit $w$ by proving   for $w^n$  a bound   with respect the metrics \be\label{MM} C(\vep,T;J^{1,2}(\OO))\cap L^2(\vep,T;W^{2,2}(\OO))\mbox{ and }w_t,\n \pi_w\in L^2(\vep,T;L^2(\OO)),\ee for arbitrary $T>\vep>0$ and then uniform in $n\in\N$ such that $\vep-\frac1n>\frac\vep2$. In particular the limit $w$ satisfies the inequality
\be\label{GWI}\ba{l}\displ \displ(t-\vep)\dm \n w(t)\dm_2^2+\mbox{$\frac12$}\intll\vep t(\tau-\vep)(\dm P\Delta w(\tau)\dm_2^2+\dm w_\tau(\tau)\dm_2^2)d\tau \VS \hskip 6cm  \leq \dm v_\circ\dm_2^2\exp[t-\vep+c\intll\vep t\dm v(\tau)\dm_\frac{2\sigma}{\sigma-2}^\frac{4\sigma}{6-\sigma}d\tau] \,.\ea \ee By virtue of estimate  \rf{EEW}, assumption for $v$, via Lemma\,\ref{EEI} we realize that $w^n\in C(0,T;J^2(\OO))$ for all $n\in\N$ with \be\label{EEW-I}\dm w^n(t)\dm_2^2+2\intll st\dm \n w^n(\tau)\dm_2^2d\tau=\dm w^n(s)\dm_2^2\leq \dm v_\circ\dm_2^2,\mbox{ for all }t>s\geq0 . \ee Now, we look for   estimates for the derivatives  of $w^n(t)$. We multiply   equation \rf{LNSM} by $P\Delta w^n$. Integrating on $(0,T)\times\OO$, we get
$$\mbox{$\frac12\frac d{dt}$}\dm\n w^n(t)\dm_2^2+\dm P\Delta w^n(t)\dm_2^2= (J_n[v](t)\cdot\n w^n(t), P\Delta w^n(t)),\mbox{ for all }t>0.$$ Analogously multiplying \rf{LNSM} by $w^n_t$ and integrating on $(0,T)\times\OO$, we get
$$\mbox{$\frac12\frac d{dt}$}\dm\n w^n(t)\dm_2^2+\dm   w^n_t(t)\dm_2^2= -(J_n[v](t)\cdot\n w^n(t),   w^n_t(t)),\mbox{ for all }t>0.$$  Applying estimate \rf{TRI}  with $q$ such that $\frac{2q}{q-2}=\sigma$ on the right hand side of the above relations with $u=P\Delta w^n$ in the first relation and with $u=w^n_t$ in the second relation, summing we get
$$\mbox{$\frac d{dt}$}\dm\n w^n(t)\dm_2^2+ \mbox{$\frac 12$}\dm P\Delta w^n(t)\dm_2^2+\mbox{$\frac 32$}\dm w^n_t(t)\dm_2^2\leq c(\dm J_n[v]\dm_\sigma^2+\dm J_n[v]\dm_\sigma^\rho)\dm \n w^n(t)\dm_2^2,$$ where, via the assumption on $v$, we have taken $\frac 3\sigma+\frac2\rho=1$ into account. Multiplying by $\tau-\vep$ and integrating over $(\vep,t)$, we deduce the boundness \be\label{FE}\ba{l}\displ(t-\vep)\dm \n w^n(t)\dm_2^2+\mbox{$\frac12$}\intll\vep t(\tau-\vep)\left[\dm P\Delta w^n(\tau)\dm_2^2+\dm w^n_\tau(\tau)\dm_2^2\right]d\tau \VS\hskip4cm\leq \exp[t-\vep+c\intll\vep t\dm J_n[v]\dm_\sigma^\rho d\tau]\intll\vep t\dm \n w^n(\tau)\dm_2^2d\tau\,.\ea\ee
By virtue of energy relation \rf{EEW-I}, and of the properties of the mollifier, estimate \rf{FE} is true for all $\vep>0$ and uniformly in $n\in\N$: 
\be\label{FE-I}\ba{l}\displ \displ(t-\vep)\dm \n w^n(t)\dm_2^2+\mbox{$\frac12$}\intll\vep t(\tau-\vep)\left[\dm P\Delta w^n(\tau)\dm_2^2+\dm w^n_\tau(\tau)\dm_2^2\right]d\tau \VS \hskip 6cm  \leq \dm v_\circ\dm_2^2\exp[t-\vep+c\intll{\vep-\frac1n} t\dm v(\tau)\dm_\frac{2\sigma}{\sigma-2}^\frac{4\sigma}{6-\sigma}d\tau] \,.\ea \ee
Estimates \rf{EEW-I} and \rf{FE-I} allow to deduce the existence of a limit $w$ belonging to \rf{MM}.   Moreover, the limit $w$ satisfies the energy equality
$$\dm w(t)\dm_2^2+2\intll st\dm \n w(\tau)\dm_2^2d\tau=\dm w(s)\dm_2^2,\mbox{ for all }t>s\geq\vep>0$$ and the energy inequality      $$\dm w(t)\dm_2^2+2\intll 0t\dm \n w(\tau)\dm_2^2d\tau\leq\dm v_\circ\dm_2^2,\mbox{ for all }t>0.$$ Finally, from estimate \rf{FE-I} for all $t>\vep$ it follows $$\displ \displ(t-\vep)\dm \n w(t)\dm_2^2+\mbox{$\frac12$}\intll\vep t(\tau-\vep)\left[\dm P\Delta w(\tau)\dm_2^2+\dm w_\tau(\tau)\dm_2^2\right]d\tau    \leq \dm v_\circ\dm_2^2\exp[t-\vep+c\intll\vep t\dm v(\tau)\dm_\frac{2\sigma}{\sigma-2}^\frac{4\sigma}{6-\sigma}d\tau] \,. $$ On the other hand, from integral equation iii) of Definition\,\ref{LNSWS} we deduce that for all $T>0$ $(w(t),\psi)\in C([0,T])$. Hence by Lemma\,\ref{EEI} we have that the limit $w$ also satisfies \rf{EE}.
In our hypotheses on $v$ weak solution to problem \rf{NS}, in accord with iii) of Definition\,\ref{LNSWS},   we  can regard it as   a   solution to \rf{LNS} with initial data $v_\circ$ and coefficient $a_2\equiv v$. So that  for problem \rf{LNS}  with $a_2=v$ we have found  two   solutions   corresponding to $v_\circ$:   $w$ and $v$. By virtue of Corollary\,\ref{C-I} they coincide\,\footnote{\; We recall that on any finite interval $(\vep,T)$ the coefficient $a_2$ in particular enjoys the property of $a_1$.}.  The proof of the lemma is completed.\chiu 
\begin{lemma}\label{LLCR}{\sl In the hypotheses of Theorem\,{\rm\ref{CT}} for all $\vep>0$   we get $\n^\ell D_t^k v\in C(\vep,T;L^2(\OO))$. In particular we get
\be\label{LLCR-I}(t-\vep) \dm v_t\dm_2\leq c\exp\big[c\intll\vep t\dm v(\tau)\dm_\sigma^\rho d\tau\big]A(\dm v_\circ\dm_2,\vep,t)
\leq B(\dm v_\circ\dm_2,\vep,t),\mbox{ for all }t>\vep\,.\ee}\end{lemma}
\Pr By virtue of \rf{LCR-I}$_1$, Theorem\,\ref{TR} and the uniqueness of strong regular solutions, we can claim that for all $k\in\N,\,\ell=1,2$ and $\vep>0$ there exists a $T_\vep$ such that $\nabla^\ell D^k_tv\in C(\vep,T_\vep;L^2(\OO))$. We can iterate the above procedure proving that   $\n^\ell D_tv\in C(\vep,T;L^2(\OO))$ for all $T>0$.  As matter of   course this procedure furnishes an estimate of $\dm v_t(t)\dm_2$ which is not uniform.  Hence in order to prove \rf{LLCR-I} we  differentiate the equation \rf{NS}$_1$ with respect to $t$, then, multiplying by $(t-\vep)^2 v_t$ and integrating on $ \OO$, we get
$$\mbox{$\frac d{dt}$}\Big[(t-\vep)^2\dm v_t(t)\dm_2^2\Big]+2(t-\vep)^2\dm\n v_t\dm_2^2=-2(t-\vep)^2(v_t\cdot\n v_t,v)+2(t-\vep)\dm v_t\dm_2^2.$$
Applying estimate \rf{TRII}, we deduce the estimate $$|(v_t\cdot\n v_t,v)|\leq c\dm v_t\dm_2^2 \dm v\dm_\sigma^\frac{2\sigma}{\sigma-3}+\mbox{$\frac14$}\dm\n v_t\dm_2^2\,.$$ Hence applying this inequality on the right hand side of the above differential equation, and integrating on $(\vep,t)$ we get
$$(t-\vep)^2\dm v_t(t)\dm_2^2\leq c\exp\big[c\intll\vep t\dm v(\tau)\dm_\sigma^\rho d\tau\big]\intll\vep t(\tau-\vep)\dm v_\tau(\tau)\dm_2^2d\tau,\mbox{ for all }t>\vep\,.$$ Taking into account \rf{EGP}, we complete the proof.  \chiu
\par {\it Proof of Theorem\,{\rm\ref{CT}}.}
\par\noindent The proof of Theorem\,\ref{CT} is an immediata consequence of the above lemmas and of the Sobolev embedding theorem. Actually, by virtue of the above lemmas and the local strong regularity Theorem\,\ref{TR}, from equation \rf{NS}$_1$ for all $\vep>0$ and  $t\in(\vep,T_\vep)$ we have  $$\dm P\Delta v(t)\dm_2 \leq \dm P(v_t+v\cdot\n v)\dm_2\leq \dm v_t\dm_2+\dm v\dm_\infty\dm\n v\dm_2\leq\dm v_t(t)\dm_2+c\dm D^2v\dm_2^\frac12\dm v\dm_6^\frac12\dm \n v\dm_2,$$ where estimating $\dm v\dm_\infty$ we have applied \rf{GNI}$_2$. By virtue of estimates \rf{EGP} and \rf{LLCR-I} the right hand side is independent of $T_\vep$, and since $v\in C((\vep,T_\vep); J^{1,2}(\OO)$ for all $\vep>0$, we can iterate  the last arguments, hence   for all $\vep>0$ and  $t>0$ we deduce  $$\dm P\Delta v(t)\dm_2 \leq \dm P(v_t+v\cdot\n v)\dm_2\leq \dm v_t\dm_2+\dm v\cdot\n v\dm_2\leq\dm v_t(t)\dm_2+c\dm D^2v\dm_2^\frac12\dm v\dm_6^\frac12\dm \n v\dm_2.$$Hence applying Lemma\,\ref{TR} and the Sobolev inequality, for all $\vep>0$ we get  
$$\dm P\Delta v\dm_2\leq \dm v_t\dm_2+c(\dm P\Delta v\dm_2+\dm \n v\dm_2)^\frac12 \dm \n v\dm_2^\frac32,\mbox{ for all }t>\vep,$$ and by the Cauchy inequality we have for all $\vep>0$ 
$$\dm P\Delta v\dm_2\leq 2\dm v_t\dm_2+c\dm \n v\dm_2^2(\dm \n v\dm_2 +1)\,,\mbox{ for all }t>\vep.$$ Finally, via \rf{EGP} and \rf{LLCR-I} for all $\vep>0$ we have 
$$(t-\vep)\dm P\Delta v\dm_2\leq cA^2(\dm v_\circ\dm_2,\vep,t)\big[(t-\vep)^{-\frac12}A(\dm v_\circ\dm_2,\vep,t)+1\big]+B(\dm v_\circ\dm_2,\vep,t),\mbox{ for all } t>\vep\,.$$ Now it is immediate to deduce the thesis.\chiu
\par{\small\bf Acknowledgements }-{\small This research was partly supported by GNFM-INdAM, and by MIUR via the PRIN 2016 ``Non-linear Hyperbolic Partial Differential Equations, Dispersive and Transport Equations: Theoretical and Applicative Aspects''.}
\end{document}